\documentclass{amsart}
\usepackage{hyperref}

\newtheorem{thm}{Theorem}
\newtheorem{lem}[thm]{Lemma}
\newtheorem{prop}[thm]{Proposition}
\theoremstyle{definition}
\newtheorem{defi}[thm]{Definition}
\newcommand{\Z}{\mathbb{Z}}

\begin{document}

\title{An algorithm for enumerating difference sets}

\author{Dylan Peifer}
\address{Department of Mathematics\\
Cornell University\\
Ithaca, NY 14853}
\email{djp282@cornell.edu}
\urladdr{pi.math.cornell.edu/~djp282}


\begin{abstract}
The {\tt DifSets} package for {\tt GAP} implements an algorithm for enumerating all difference sets in a group up to equivalence and provides access to a library of results.
The algorithm functions by finding difference sums, which are potential images of difference sets in quotient groups of the original group, and searching their preimages.
In this way, the search space can be dramatically decreased, and searches of groups of relatively large order (such as order 64 or order 96) can be completed.
\end{abstract}

\maketitle

\section{Introduction}\label{sec:introduction}

Let $G$ be a finite group of order $v$ and $D$ a subset of $G$ with $k$ elements.
Then $D$ is a \emph{$(v, k, \lambda)$-difference set} if each nonidentity element of $G$ can be written as $d_id_j^{-1}$ for $d_i, d_j \in D$ in exactly $\lambda$ different ways.
Difference sets were first studied in relation to finite geometries \cite{Singer1938} and have connections to symmetric designs, coding theory, and many other fields of mathematics \cite{MoorePollatsek2013, DavisJedwab1996, CRC, BeJuLe1999}.

Large libraries of difference sets are useful for developing conjectures and building examples.
Gordon \cite{LaJolla} provides an extensive library of difference sets in abelian groups, but has no results for non-abelian groups, which do show distinct behavior \cite{Smith1995}.
A wide variety of techniques can be used to construct difference sets for these libraries (see, for example, \cite{Dillon1985} and \cite{DavisJedwab1997}), but fully enumerating all difference sets in a given group requires some amount of exhaustive search, which can quickly become computationally infeasible.
Kibler \cite{Kibler1978} performed the first major exhaustive enumeration of difference sets, and considered groups where difference sets could be found with $k < 20$.
In recent years, AbuGhneim \cite{AbuGhneim2013, AbuGhneim2016} performed almost complete enumerations for all groups of order 64, and several authors have found difference sets in groups of order 96 \cite{GoMaVu2005, GoMaVu2007, AbuGhneimSmith2007}.

The {\tt DifSets} package for the {\tt GAP} \cite{GAP4} computer algebra system efficiently and generally implements the techniques used by these and many other authors to exhaustively enumerate all difference sets up to equivalence in a group.
With the package loaded, a search of a given group can be performed with a single command.

\begin{verbatim}
gap> DifferenceSets(CyclicGroup(7));
[ [ 1, 2, 4 ] ]
\end{verbatim}

The package has been used to give the first complete enumeration of all difference sets up to equivalence in groups of order 64 and 96, and in total provides a library of results for 1006 of the 1032 groups of order less than 100.
Results are organized by their ID in the {\tt SmallGroups} \cite{SmallGrp} library and can be easily loaded by {\tt GAP}.

\begin{verbatim}
gap> LoadDifferenceSets(16, 5);  # results for SmallGroup(16, 5)
[ [ 1, 2, 3, 4, 8, 15 ], [ 1, 2, 3, 4, 11, 13 ] ]
\end{verbatim}

The ease of use of these top-level functions is the primary interface difference between the {\tt DifSets} package and the {\tt RDS} \cite{RDS} package, a similar {\tt GAP} package that provides a variety of tools to search for difference sets.
The functions involving coset signatures in {\tt RDS} provide similar functionality to the {\tt DifSets} package, but require substantial user interaction to perform efficient searches, and are not feasible for searching most groups of order 64 and 96.
In addition, {\tt RDS} provides no precomputed results, though it does provide significant additional functionality related to relative difference sets, partial difference sets, and projective planes.

\section{Difference Sets}\label{sec:difsets}

For notational purposes it is useful to consider a subset $D \subseteq G$ as an element of the group ring $\Z[G]$.
We will abuse notation to define the group ring elements
\[
G = \sum_{g \in G} g, \quad D = \sum_{d \in D} d, \quad D^{(-1)} = \sum_{d \in D} d^{-1}, \quad gD = \sum_{d \in D} gd, \quad D^\phi = \sum_{d \in D} \phi(d)
\]
where $g \in G$ and $\phi$ is a homomorphism from $G$.
Then the statement that $D$ is a $(v, k, \lambda)$-difference set is equivalent to the equation
\[
DD^{(-1)} = (k-\lambda)1_{G} + \lambda G
\]
where $D$ is an element of $\Z[G]$ with coefficients in $\{0, 1\}$.
With this definition it is a quick exercise to prove the following (see page 298 of \cite{BeJuLe1999} and Theorem 4.2 and 4.11 of \cite{MoorePollatsek2013}).

\begin{prop}\label{prop:equivsets}
Let $G$ be a group of order $v$.
Then
\begin{enumerate}
\item Any one element subset of $G$ is a $(v, 1, 0)$-difference set.
\item The complement of a $(v, k, \lambda)$-difference set in $G$ is a $(v, v-k, \lambda+v-2k)$-difference set in $G$.
\item If $D$ is a $(v, k, \lambda)$-difference set in $G$, $g \in G$, and $\phi \in \mathrm{Aut}(G)$, then $gD^\phi$ is also a $(v, k, \lambda)$-difference set in $G$.
\end{enumerate}
\end{prop}

In addition, an immediate consequence of the definition is that $k(k-1) = \lambda (v-1)$ for any valid set of parameters of a difference set, so that for a given value of $v$ there are typically only a few possible values of $k$ and $\lambda$.
More sophisticated results, such as the Bruck-Ryser-Chowla theorem, can reduce the number of possibilities even further.

As a result of Proposition~\ref{prop:equivsets}, in enumerating difference sets we ignore the trivial one element difference sets, only take the smaller of each complementary pair of sets, and only consider sets distinct up to an equivalence given by part (3).

\begin{defi}\label{defi:equivsets}
Let $D_1$ and $D_2$ be difference sets in $G$.
Then $D_1$ and $D_2$ are \emph{equivalent difference sets} if $D_1 = gD_2^\phi$ for some $g \in G$ and $\phi \in \mathrm{Aut}(G)$.
\end{defi}

In the {\tt DifSets} package, difference sets are stored as lists of integers.
These integers represent indices in the list returned by the {\tt GAP} function {\tt Elements(G)}, which is a sorted\footnote{Element comparison (and thus the list {\tt Elements(G)}) is instance-independent in {\tt GAP} for permutation and pc groups, which includes, for example, all groups in the {\tt SmallGroups} library.} list of elements of the group {\tt G}.
For example, consider the group $C_7 = \langle x | x^7 = 1 \rangle$.
In {\tt GAP} we have
\begin{verbatim}
gap> C7 := CyclicGroup(7);;
gap> Elements(C7);
[ <identity> of ..., f1, f1^2, f1^3, f1^4, f1^5, f1^6 ]
\end{verbatim}
where clearly {\tt f1} is the generator corresponding to our $x$.
Then the subset $D = \{x, x^2, x^4\}$ corresponds to the set consisting of the second, third, and fifth elements of {\tt Elements(C7)}, which we can represent in indices as {\tt [2, 3, 5]}.
We can check that this is a difference set and also note that it is equivalent to the difference set $xD = \{x^2, x^3, x^5\}$, which is represented as {\tt [3, 4, 6]}.
\begin{verbatim}
gap> IsDifferenceSet(C7, [2, 3, 5]);
true
gap> IsEquivalentDifferenceSet(C7, [2, 3, 5], [3, 4, 6]);
true
\end{verbatim}

\section{Difference Sums}

A basic method for enumerating all difference sets in a group $G$ is to enumerate all subsets of $G$ and check if each is a difference set by definition.
But since the number of subsets in a group is exponential in its order, we cannot feasibly enumerate and test all subsets for groups of even a modest size.
The key to decreasing the search space is the following well-known lemma, which motivates our definition of a \emph{difference sum}\footnote{Concepts similar to difference sums are elsewhere referred to as difference lists, intersection numbers, or signatures.
However, difference sums require both a group $G$ and normal subgroup $N$, not just the group structure of the quotient $G/N$ used in some other definitions.
This precision is needed for specifying induced automorphisms in Definition~\ref{defi:equivsums} so that we can prove Lemma~\ref{lem:equivsums}.}.

\begin{lem}\label{lem:image}
Suppose $D$ is a $(v, k, \lambda)$-difference set in $G$ and $\theta$ is a homomorphism of $G$ with $|\mathrm{ker}(\theta)| = w$. Let $S = D^\theta$ and $H = G^\theta$.
Then
\[
SS^{(-1)} = (k - \lambda)1_H + \lambda w H.
\]
\end{lem}

\begin{defi}\label{defi:difsum}
Given a finite group $G$ and normal subgroup $N$, a $(v, k, \lambda)$-\emph{difference sum} is an element $S$ of $\Z[G/N]$ such that $SS^{(-1)} = (k - \lambda)1_{G/N} + \lambda |N| G/N$ and the coefficients of $S$ have values in $\{0, 1, \dots, |N|\}$.
\end{defi}

By construction, any difference set in $G$ induces difference sums under the natural projection in quotients of $G$, as seen in Figure~\ref{fig:induce}.
Precisely, we have

\begin{figure}
\renewcommand{\arraystretch}{1.3}
\begin{tabular}{|c|c|c|c|c|c|c|c|c|c|c|c|c|c|c|cc}
\cline{1-15}
1 & 1 & 1 & 0 & 0 & 1 & 0 & 1 & 0 & 0 & 1 & 0 & 0 & 0 & 1 & \quad & $G/N_3$ \\
\cline{1-15}
\multicolumn{3}{|c|}{3} & \multicolumn{3}{|c|}{1} & \multicolumn{3}{|c|}{1} & \multicolumn{3}{|c|}{1} & \multicolumn{3}{|c|}{1} & \quad & $G/N_2$ \\
\cline{1-15}
\multicolumn{15}{|c|}{7} & \quad & $G/N_1$ \\
\cline{1-15}
\end{tabular}
\caption{A difference set of size 7 in the group $G$ of order 15 and the difference sums it induces in $G/N_i$ where $G = N_1 \triangleright N_2 \triangleright N_3 = \{1\}$.
Each row in the diagram is a group, with each block a coset.
\label{fig:induce}}
\end{figure}

\begin{lem}\label{lem:setinduce}
Suppose $G$ is a finite group with normal subgroup $N$ and natural projection $\pi : G \to G/N$.
Then any $(v, k, \lambda)$-difference set $D$ in $G$ induces a $(v, k, \lambda)$-difference sum $D^\pi$ in $G/N$.
\end{lem}

\begin{lem}\label{lem:suminduce}
Suppose $G$ is a finite group with normal subgroups $N_1$ and $N_2$ such that $N_2 \subseteq N_1$ and $\pi : G/N_2 \to G/N_1$ is the natural projection.
Then any $(v, k, \lambda)$-difference sum $S$ in $G/N_2$ induces a $(v, k, \lambda)$-difference sum $S^\pi$ in $G/N_1$.
\end{lem}

Lemma~\ref{lem:setinduce} means that our search for difference sets only requires checking the subsets of $G$ that induce difference sums in some quotient.
In finding these difference sums, Lemma~\ref{lem:suminduce} additionally allows us to only test sums that induce difference sums in further quotients.
In each case the search space is dramatically decreased.
Since our search is for difference sets up to equivalence, we also define a complementary equivalence of difference sums such that equivalent difference sums are induced by equivalent collections of difference sets.

\begin{defi}\label{defi:equivsums}
Let $S_1$ and $S_2$ be difference sums in $G/N$.
Then $S_1$ and $S_2$ are \emph{equivalent difference sums} if $S_1 = gS_2^\phi$ for some $g \in G/N$ and $\phi$ an automorphism of $G/N$ induced by an automorphism of $G$.
\end{defi}

\begin{lem}\label{lem:equivsums}
Suppose $S_1$ and $S_2$ are equivalent difference sums in $G/N$.
Then if $D_1$ is any difference set in $G$ that induces $S_1$, there exists a difference set $D_2$ in $G$ that induces $S_2$ such that $D_1$ and $D_2$ are equivalent.
\end{lem}

In the {\tt DifSets} package, difference sums are stored as lists of integers representing the coefficients of the group ring elements, with position in the list given by the position of the coset in the list returned by the {\tt GAP} function {\tt Elements(G/N)}.
For example, {\tt [3, 1, 1, 1, 1]} represents a difference sum in {\tt SmallGroup(15, 1)} mod its normal subgroup of order 3 with coefficient 3 on the identity coset and coefficient 1 on all other cosets.

\begin{verbatim}
gap> G := SmallGroup(15, 1);; N := NormalSubgroups(G)[2];;
gap> IsDifferenceSum(G, N, [3, 1, 1, 1, 1]);
true
\end{verbatim}

\section{Algorithm}\label{sec:algorithm}

The basic structure of the algorithm is to start at the bottom of Figure~\ref{fig:induce} and travel upwards.
Given a group $G$, first compute $v = |G|$ and then find all values of $k$ that give solutions satisfying the Bruck-Ryser-Chowla theorem to the equation $k(k-1) = \lambda (v-1)$ mentioned in Section~\ref{sec:difsets}.
For example,

\begin{verbatim}
gap> G := SmallGroup(15, 1);;
gap> PossibleDifferenceSetSizes(G);
[ 7 ]
\end{verbatim}

Each value of $k$ will be handled separately.
The algorithm starts with the normal subgroup $N_1 = G$, where the only difference sum of size $k$ in $G/N_1 = \{1\}$ is {\tt [k]}.

\begin{verbatim}
gap> N1 := G;;
gap> difsums := [ [7] ];;
\end{verbatim}

Given a normal subgroup $N_2$ of $G$ such that $N_2 \subseteq N_1$, first enumerate all preimages in $G/N_2$ of current difference sums in $G/N_1$ and return those that are themselves difference sums.
Then remove all but one representative of each equivalence class from this collection.

\begin{verbatim}
gap> N2 := NormalSubgroups(G)[2];;
gap> difsums := AllRefinedDifferenceSums(G, N1, N2, difsums);
[ [ 1, 1, 1, 1, 3 ], [ 1, 1, 1, 3, 1 ], [ 1, 1, 3, 1, 1 ],
  [ 1, 3, 1, 1, 1 ], [ 3, 1, 1, 1, 1 ] ]
gap> difsums := EquivalentFreeListOfDifferenceSums(G, N2, difsums);
[ [ 3, 1, 1, 1, 1 ] ]
\end{verbatim}

In the general case, the above step is repeated along a chief series $G = N_1 \triangleright \dots \triangleright N_r = \{1\}$ of $G$ with $N_{r-1}$ a nontrivial normal subgroup of minimal possible size in $G$.
At $N_{r-1}$, enumerate sets and remove equivalents to leave the final result.

\begin{verbatim}
gap> difsets := AllRefinedDifferenceSets(G, N2, difsums);
[ [ 1, 2, 4, 3, 8, 11, 12 ], [ 1, 2, 4, 3, 10, 13, 12 ], 
  [ 1, 2, 4, 5, 6, 9, 14 ], [ 1, 2, 4, 5, 10, 13, 14 ], 
  [ 1, 2, 4, 7, 6, 9, 15 ], [ 1, 2, 4, 7, 8, 11, 15 ] ]
gap> difsets := EquivalentFreeListOfDifferenceSets(G, difsets);
[ [ 1, 2, 4, 7, 8, 11, 15 ] ]
\end{verbatim}

These steps are encapsulated in the function {\tt DifferenceSets} mentioned in Section~\ref{sec:introduction}, with two modifications.
First, since every difference set is equivalent to some difference set containing the identity, the algorithm does not enumerate some preimages that are guaranteed to be equivalent to others.
Second, the final elimination of all but one representative of equivalence classes of difference sets uses the {\tt SmallestImageSet} function \cite{Linton2004} from the {\tt GAP} package {\tt GRAPE} \cite{GRAPE}.
Although roughly 20\% slower than the function given above for most cases, {\tt SmallestImageSet} gives a unique minimal result and handles groups with large automorphism groups much more efficiently.

\section{Results}

The {\tt DifSets} package successfully computed results for 1006 of the 1032 groups of order less than 100, including all groups of order 64 and 96.
Full results with timings and comments can be found in the package and its documentation.
Here we include a summary for order 64 and 96.
All computations were performed with {\tt GAP} 4.9.1 on a 4.00GHz i7-6700K using 8GB of RAM.

\begin{center}
\begin{tabular}{ccccc}
Order & Groups & Difference Sets & Median Time per Group & Total Time \\
\hline
64 & 267 & 330159 & 0.415 hours & 295.811 hours \\
96 & 231 & 2627 & 3.133 hours & 1568.746 hours
\end{tabular}
\end{center}

Timing comparisons with the {\tt RDS} package mentioned in Section~\ref{sec:introduction} are difficult since {\tt RDS} provides a variety of tools rather than a single algorithm.
Ordered coset signatures in {\tt RDS} correspond to difference sums in {\tt DifSets}, but, unlike difference sums, coset signatures cannot be refined through multiple stages, which makes the generation of good coset signatures in {\tt RDS} infeasible for most order 64 and order 96 groups.
However, if an ordered signature is available, building difference sets through partial difference sets in {\tt RDS} can in some cases be much faster than searching the corresponding difference sum using {\tt DifSets}.
In particular, replacing the final step in Section~\ref{sec:algorithm} with a search using {\tt RDS} can significantly improve times for some groups of order 96.
Further work to combine the refining of difference sums used by {\tt DifSets} and the generation of difference sets through partial difference sets used by {\tt RDS} could lead to significantly better times than either package could manage alone.

\section*{Acknowledgements}

The author thanks Ken Smith, Alexander Hulpke, and an anonymous reviewer for helpful comments that improved the {\tt DifSets} package and this article.

\bibliographystyle{plain}
\bibliography{Peifer2018}

\end{document}